\documentclass{amsart}
\usepackage{amssymb,euscript,amsmath, mathrsfs}
\usepackage[dvips]{graphicx}
\usepackage[dvips]{color}

\newcounter{ENUM}

\input xy
\xyoption{all}
\CompileMatrices

\def\<{\langle}
\def\>{\rangle}
\def\0{{{\bf 0}}}

\def\OO{{\mathcal O}}

\def\QQ{{\mathbb Q}}

\def\TT{{\mathbb T}}

\def\rhobar{{\overline{\rho}}}
\def\rhobar{{\overline{\rho}}}
\def\pibar{{\overline{\pi}}}

\def\abs{\mid~\;~\mid}
\def\unif{\varpi}

\newcommand{\supp}{\operatorname{Supp}}

\newcommand{\St}{\operatorname{St}}

\newcommand{\env}{\operatorname{env}}

\newcommand{\sss}{\mbox{\rm \tiny ss}}

\newcommand{\et}{\mbox{\rm \tiny \'et}}

\def\ZZ{{\mathbb Z}}

\def\gen{\operatorname{gen}}

\def\Ext{\operatorname{Ext}}

\def\Hom{\operatorname{Hom}}

\def\End{\operatorname{End}}

\def\GL{\operatorname{GL}}

\newcommand{\margh}[1]{}

\newtheorem{thm}{Theorem}[section]
\newtheorem{prop}[thm]{Proposition}
\newtheorem{lemma}[thm]{Lemma}

\theoremstyle{definition}

\numberwithin{equation}{section}

\begin{document}
\title{On the modified mod $p$ local Langlands correspondence for $\GL_2(\QQ_{\ell})$}
\author{David Helm}
\subjclass[2000]{11F70 (Primary), 22E50 (Secondary)}

\maketitle

We give an explicit description of the modified mod $p$ local Langlands correspondence
for $\GL_2(\QQ_{\ell})$ of~\cite{emerton-helm}, Theorem 5.1.5, where $p$ is an odd
prime different from $\ell$.

\section{Introduction}
In~\cite{emerton-helm}, Matthew Emerton and the author introduce a
``modified mod $p$ local Langlands correspondence,'' a ``mod $p$''
version of the local Langlands correspondence that is well-behaved under
specialization and has useful applications to the cohomology of modular curves
and the ``local Langlands correpondence in families'' of~\cite{emerton-helm}.
Section 5 of~\cite{emerton-helm} gives a general characterization of this mod $p$
correspondence in terms of its basic properties.  If one restricts to the group $\GL_2$,
it is easy in most cases to go from
this list of characterizing properties to an explicit description of this correspondence.
These easy cases are discussed in detail in section 5.2 of~\cite{emerton-helm}.  When $p$ is odd
the cases discussed come close to a complete description of the correspondence, but omit
certain more difficult special cases.  The purpose of this note is to explicitly describe the
correspondence in these more difficult cases and thus complete the description of the modified mod $p$
local Langlands correspondence for $\GL_2$ and odd $p$. 

Let $F$ be a finite extension of $\QQ_{\ell}$ whose residue field has order $q$, 
let $p$ be an odd prime distinct from $\ell$,
and let $k$ be a finite field of characteristic $p$.
The modified mod $p$ local Langlands correspondence is an association $\rhobar \mapsto \pibar(\rhobar)$,
where $\rhobar: G_F \rightarrow \GL_n(k)$ is a continuous $n$-dimensional representation of the absolute Galois
group of $F$, and $\pibar(\rhobar)$ is a finite length indecomposable smooth representation of
$\GL_n(F)$.  Its interest arises from its nice behaviour under specialization, which
we discuss below, and also from the fact that it arises ``in nature'' in the cohomology
of the tower of modular curves.  Indeed, in~\cite{emerton-lg}, Emerton considers the following situation:

Let $\Sigma$ be a finite set of primes containing $p$, and let $H^1_{\Sigma}$ be the direct
limit:
$$\lim_{\supp N \subseteq \Sigma} H^1_{\et}(X(N)_{\overline{\QQ}}, k).$$
(Here $\supp N$ denotes the set of primes dividing an integer $N$, so the limit is
over $N$ divisible only by primes in $\Sigma$, ordered by divisibility.)
The space $H^1_{\Sigma}$ acquires actions of $G_{\QQ}$, of $\GL_2(\QQ_p)$, and of $\GL_2(\QQ_{\ell})$
for $\ell \neq p$, as well as of the Hecke operators $T_r$ for $r \notin \Sigma$ and the diamond operators
$\<d\>$ for $d$ not divisible by any prime of $\Sigma$.  Let $\TT_{\Sigma}$ be the subalgebra
of $\End_k(H^1_{\Sigma})$ generated by these Hecke operators and diamond operators.

Let $\rhobar: G_{\QQ} \rightarrow \GL_2(k)$ be a modular Galois representation unramified outside $\Sigma$.
Then there is a maximal ideal ${\mathfrak m}$ of $\TT_{\Sigma}$ attached to $\rhobar$, and
(under certain hypotheses on the local behavior of $\rhobar$ at $p$),
Emerton has shown (\cite{emerton-lg}, Theorem 6.2.13 and Proposition 6.1.20)
that $H^1_{\sigma}[{\mathfrak m}]$ is a product of ``local factors'':
$$H^1_{\sigma}[{\mathfrak m}] \cong \rhobar \otimes \pi_p \otimes
\bigotimes_{\ell \neq p, \ell \in \Sigma} \pi_{\ell}$$
where $\pi_p$ is attached to $\rhobar|_{G_{\QQ_p)}}$ by considerations from the $p$-adic Langlands
program (see~\cite{emerton-lg}, section 3, for details)
and each $\pi_{\ell}$ is the representation $\pibar(\rhobar|_{G_{\QQ_{\ell}}})$ attached to
the restriction of $\rhobar$ to a decomposition group at $\ell$ via the modified mod $p$ local
Langlands correspondence for $\GL_2(\QQ_{\ell})$.  Thus an explicit description of the modified mod
$p$ local Langlands correspondence for $\GL_2(\QQ_{\ell})$ gives an explicit description
of the action of $\GL_2(\QQ_{\ell})$ on the cohomology of the modular tower.

We now recall more precisely the definition of the modified mod $p$ local Langlands
correspondence.  The starting point is the characteristic zero ``generic local Langlands
correspondence'' of Breuil-Schneider~\cite{BS}.  We refer the reader to sections 4.2 and 4.3
of~\cite{emerton-helm} for the basic properties of this correspondence.  In particular, this
correspondence associates to any $n$-dimensional
Frobenius-semisimple Weil-Deligne representation $(\rho,N)$
over a field $K$ containing $\QQ_p$ an indecomposable (but often reducible) admissible
representation $\pi(\rho,N)$ of $\GL_n(F)$.

The representations $\pi(\rho,N)$ have several nice properties.  In particular, they are
{\em essentially AIG} representations, a concept introduced in section 3.2 of~\cite{emerton-helm}.
A smooth representation $\pi$ of $\GL_n(F)$ over a field $K$ is called essentially AIG if
\begin{itemize}
\item the socle of $\pi$ is absolutely irreducible and generic,
\item $\pi$ contains no generic irreducible subquotients other than its socle, and
\item $\pi$ is the sum of its finite length submodules.
\end{itemize}

Such representations have several useful properties.  In particular, their only endomorphisms are scalars
(\cite{emerton-helm}, Lemma 3.2.3), any submodule of an essentially AIG representation is essentially AIG,
any nonzero map of essentially AIG representations is an embedding, and such an embedding, if it
exists at all, is unique up to a scalar factor (\cite{emerton-helm}, Lemma 3.2.2).
Moreover, if $\pi$ is an absolutely irreducible generic representation of $\GL_n(F)$, then there is an
essentially AIG representation $\env(\pi)$, known as the essentially AIG envelope of $\pi$, such that
the socle of $\env(\pi)$ is isomorphic to $\pi$ and any essentially AIG representation $\pi'$ with socle
isomorphic to $\pi$ embeds in $\env(\pi)$ (\cite{emerton-helm}, Proposition 3.2.7).  Moreover, 
all the subquotients of $\env(\pi)$ (or, more generally, of any essentially AIG representation) have the
same supercuspidal support (\cite{emerton-helm}, Corollary 3.2.14).

A final useful property of essentially AIG representations is that they contain distinguished lattices (up to
homothety).  In particular, let $\OO$ be a discrete valuation ring with residue field $k$ and field of fractions
$K$, and let $\pi$ be an essentially AIG representation over $K$.  Suppose further that $\pi$ is $\OO$-integral;
that is, contains a $\GL_n(F)$-invariant $\OO$-lattice.  Then there is an $\OO$-lattice $\pi^{\circ}$ in
$\pi$, unique up to homothety, such that $\pi^{\circ} \otimes_{\OO} k$ is essentially AIG (\cite{emerton-helm}, 
Proposition 3.3.2).

This last property is crucial because it allows for a definition of the modified mod $p$ local Langlands
correspondence via ``compatibility with reduction mod $p$'' from the characteristic zero correspondence of
Breuil-Schneider.  More precisely, one has:

\begin{thm}[\cite{emerton-helm}, Theorem 5.1.5]
Let $k$ be a finite field of characteristic $p$.  There is a map $\rhobar \mapsto \pibar(\rhobar)$
from the set of isomoprhism classes of continuous representations $G_F \rightarrow \GL_n(k)$ to the set
of isomorphism classes of finite length admissible smooth $\GL_n(F)$-representations over $k$, uniquely
determined by the following three conditions:
\begin{enumerate}
\item For any $\rhobar$, the associated $\GL_n(F)$-representation $\pibar(\rhobar)$ is essentially AIG.
\item If $K$ is a finite extension of $\QQ_p$, with ring of integers $\OO$ and residue field $k'$ containing $k$,
$\rho: G_F \rightarrow \GL_n(\OO)$ is a continuous representation lifting $\rhobar \otimes_k k'$, and $\pi^{\circ}$
is the unique $\OO$-lattice in $\pi(\rho)$ such that $\pi^{\circ} \otimes_{\OO} k'$ is essentially AIG,
then there is an embedding
$$\pi^{\circ} \otimes_{\OO} k' \hookrightarrow \pibar(\rhobar) \otimes_k k'.$$
\item The representation $\pibar(\rhobar)$ is minimal with respect to satisfying conditions (1) and (2);
that is, for any $\rhobar$, and any representation $\pi$ satisfying these conditions with respect to
$\rhobar$, there is an embedding of $\pibar(\rhobar)$ in $\pi.$
\end{enumerate}
\end{thm}

The construction of $\pibar(\rhobar)$ is straightforward.  One shows that for a given $\rhobar$, and any lift
$\rho$ of $\rhobar$ as in condition (2), the socle of $\pi^{\circ} \otimes_{\OO} k'$ is the unique absolutely
irreducible generic representation $\pibar^{\gen}$ of $\GL_n(F)$ whose supercuspidal support corresponds to $\rhobar^{\sss}$
under the mod $p$ semisimple local Langlands correspondence of Vigneras~\cite{vigss}.
Thus $\pi^{\circ} \otimes_{\OO} k'$ embeds in the essentially AIG envelope $\env(\pibar^{\gen})$, so that
$\env(\pibar^{\gen})$ satisfies conditions (1) and (2) of the theorem, but may be too large.  One obtains
$\pibar(\rhobar)$ by taking the sum, inside $\env(\pibar^{\gen}) \otimes_k \overline{k}$, of the subobjects
$\pi^{\circ} \otimes_{\OO} \overline{k}$ over all lifts $\rho$ as in (2), and descending from $\overline{k}$ to $k$.

When $n=2$ and $p$ is odd, this perspective is all that one needs to explicitly describe the modified
mod $p$ local Langlands correspondence.  In particular, when $\rhobar^{\sss}$ is not a twist of
the direct sum $1 \oplus \overline{\omega}$, where $\overline{\omega}$ is the mod $p$ cyclotomic character, then
$\env(\pibar^{\gen})$ is irreducible, and thus the inclusions 
$\pibar^{gen} \subseteq \pibar(\rhobar) \subseteq \env(\pibar^{\gen})$
are all equalities.  When $\rhobar^{\sss}$ is a twist of $1 \oplus \overline{\omega}$, the situation is slightly more
complicated, but still easy as long as $q$ is not congruent to $\pm 1$ modulo $p$.  We refer the reader to section 5.2
of~\cite{emerton-helm} for details.

In section 2 we cover the case when $q$ is congruent to $-1$ modulo $p$, and
$\rhobar^{\sss}$ is a twist of $1 \oplus \omega$.  This case was worked out independently by
Emerton in unpublished work.  It is similar to the case when $q$ is not congruent to $\pm 1$ modulo $p$,
but is slightly more complicated because $\env(\pibar^{gen})$ has length $3$ instead of $2$.

The remaining sections are devoted to the case when $q$ is congruent to $1$ modulo $p$, and $\rhobar^{\sss}$ is
a twist of $1 \oplus \overline{\omega}.$  (Note that $\overline{\omega}$ is trivial in this setting.)  This case is
the most difficult because in this setting
there is a one-parameter family of representations $\rhobar$ whose semisimplification is trivial.  It turns out
(Theorem~\ref{thm:main2})
in this case that the modified mod $p$ local Langlands correspondence is sensitive enough to distinguish between
these distinct extensions.  This is in stark contrast to the situation in characteristic zero, where the
Breuil-Schneider correspondence is insensitive to Frobenius-semisimplification.  This ``extra sensitivitity'' is quite
striking and would be worth investigating in situations where $n$ is greater than two.

{\em Acknowledgements}
The results in this paper grew out of a series of discussions with Matthew Emerton, and I am indebted to him
for his ideas and suggestions.  The paper was partially supported by NSF grant DMS-1161582.

\section{$q \equiv -1$ mod $p$}

In this section we write $G$ for $\GL_2(F)$, for conciseness.  Let $B$ be the standard Borel
subgroup of $G$.
Suppose that $q$ is congruent to $-1$ modulo $p$, and that
$\rhobar^{\sss}$ is a twist of $1 \oplus \overline{\omega}$.  Since the modified mod $p$
local Langlands correspondence is compatible with twisting, we may assume that $\rhobar^{\sss}$
is equal to $1 \oplus \overline{\omega}$.  The semisimple mod $p$ local Langlands correspondence
of Vigneras then shows that $\pibar^{\gen}$, and indeed every Jordan-H\"older constituent of
$\env(\pibar^{\gen})$, has supercuspidal support given by the two characters $\abs^{\pm \frac{1}{2}}$.
Thus every Jordan-H\"older constituent of $\env(\pibar^{\gen})$ is also a Jordan-Holder
constituent of the normalized parabolic induction $i_B^G \abs^{\frac{1}{2}} \otimes \abs^{-\frac{1}{2}}$.
There are three such constituents in this setting (c.f. Example II.2.5 of~\cite{vigbook}):
the trivial character of $G$, the character $\abs \circ \det$ (which has values in $\pm 1$ because
of our assumption on $q$), and a cuspidal subquotient, which is the unique generic subquotient
and is thus isomorphic to $\pibar^{\gen}$.  More precisely, we have exact sequences:
$$0 \rightarrow W \rightarrow i_B^{\GL_2(F)} \abs^{-\frac{1}{2}} \otimes \abs^{\frac{1}{2}}
\rightarrow \abs \circ \det \rightarrow 0$$
$$0 \rightarrow 1_G \rightarrow W \rightarrow \pibar^{\gen}\rightarrow 0$$
for a suitable representation $W$.  Both of these sequences are nonsplit, as $\pibar^{\gen}$
is cuspidal and thus is neither a subobject nor a quotient of any parabolic induction.

\begin{lemma} \label{lem:ext1}
Any nonsplit extension of $\pibar^{\gen}$ by the trivial character of $G$ is isomorphic to $W$.
\end{lemma}
\begin{proof}
Let $W'$ be such an extension.
The parabolic restriction $r_{\GL_2(F)}^B W'$
is isomorphic to $\abs^{-\frac{1}{2}} \otimes \abs^{\frac{1}{2}}$.  As parabolic induction
is a right adjoint to parabolic restriction, this isomorphism gives rise to a nonzero
map:
$$W' \rightarrow i_B^G \abs^{-\frac{1}{2}} \otimes \abs^{\frac{1}{2}}$$
It is clear that such a map must be injective with image $W$.
\end{proof}

Twisting by $\abs \circ \det$ we find that $W \otimes (\abs \circ \det)$ is the
unique nonsplit extension of $\pibar^{\gen}$ by $\abs \circ \det$.

As $\pibar^{\gen}$ is self-dual, it is clear that $W^{\vee}$ and $W^{\vee} \otimes (\abs \circ \det)$
are essentially AIG representations with socle $\pibar^{\gen}$.  It follows that the representation
$V$ obtained as the pushout of the diagram:
$$
\begin{array}{ccc}
W^{\vee} & \rightarrow & V\\
\uparrow & & \uparrow\\
\pibar^{\gen} & \rightarrow & W^{\vee} \otimes (\abs \circ \det)
\end{array}
$$
is also essentially AIG.  Note that $V$ is an extension of $1_G \oplus (\abs \circ \det)$
by $\pibar^{\gen}$.

\begin{prop}
The representation $V$ is an essentially AIG envelope of $\pibar^{\gen}$.
\end{prop}
\begin{proof}
We must show that $V$ is not properly contained in an essentially AIG representation
$V'$.  Suppose there were such a $V'$.  Then the socle of $V'/V$ is a contains no generic summand,
and is thus isomorphic to a direct sum of characters, each of which is either trivial
or isomorphic to $\abs \circ \det$. Observing that $V$ is isomorphic to $V \otimes (\abs \circ \det)$
we can ensure, twisting $V'$ if necessary, that $V'/V$ contains a one-dimensional subspace on
which $G$ acts trivially.  Let $V''$ be the preimage of this subspace under the surjection
$$V' \rightarrow V'/V.$$
Then $V''$ is an essentially AIG representation containing $V$, such that $V''/V$ is the
character $1_G$, and it suffices to show that such a representation cannot exist.

Note that $V''$ is essentially AIG, and hence its only endomorphisms are scalars.  In
particular the center of $G$ acts on $V''$ by a character, and this character must be trivial
since the center of $G$ acts trivially on $\pibar^{\gen}$.  On the other hand,
$V''/\pibar^{\gen}$ is en extension of $1_G$ by $1_G \oplus (\abs \circ \det)$, and it is
easy to see that such an extension must split if the center of $G$ acts trivially.

Thus $V''$ is an extension of $1_G \oplus 1_G \oplus (\abs \circ \det)$ by $\pibar^{\gen}$.
Appliying duality to Lemma~\ref{lem:ext1} shows that there is a unique nonsplit extension
of $1_G$ by $\pibar^{\gen}$, and one deduces from this that the socle of $V''$ contains at least
one copy of $1_G$, contradicting the hypothesis that $V''$ was essentially AIG.
\end{proof}

We now turn to understanding the modified mod $p$ local Langlands correspondence.  If $\rhobar^{\sss}$
is isomorphic to $1 \oplus \overline{\omega}$, then $\rhobar$ is either a nonsplit extension
of $\overline{\omega}$ by $1$, a nonsplit extension of $1$ by $\overline{\omega}$, or the direct sum
$1 \oplus \overline{\omega}$.

By contrast, let us consider the representations $\rho: G_F \rightarrow \GL_2(K)$ whose mod $p$
reduction has semisimplification $1 \oplus \overline{\omega}$.  There are several cases:
\begin{enumerate}
\item $\rho$ is irreducible, in which case its reduction modulo $p$ can be any of the three possibilities
described above.
\item $\rho$ is a nonsplit extension of $\chi_1$ by $\chi_2$, where the mod $p$ reduction of $\chi_1$
is trivial and the mod $p$ reduction of $\chi_2$ is $\overline{\omega}$.  In this case the mod $p$ reduction
of $\rho$ has a subrepresentation isomorphic to $\overline{\omega}$, and thus cannot be a nonsplit extension of 
$\overline{\omega}$ by $1$.
\item $\rho$ is a nonsplit extension of $\chi_2$ by $\chi_1$, where the mod $p$ reduction of $\chi_1$
is trivial and the mod $p$ reduction of $\chi_2$ is $\overline{\omega}$.  In this case the mod $p$
reduction of $\rho$ cannot be a nonsplit extension of $1$ by $\overline{\omega}$.
\item $\rho$ is a direct sum of two characters.  In this case the mod $p$ reduction of $\rho$ must be
the direct sum $1 \oplus \overline{\omega}$.
\end{enumerate}

It is straightforward to describe $\pi(\rho)$, and the reduction $\pi(\rho)^{\circ} \otimes_{\OO} k'$,
in each of the above cases:
\begin{enumerate}
\item $\pi(\rho)$ is irreducible and cuspidal.  In this case there is a unique homothety class
of lattices in $\pi^{\rho}$, and the reduction $\pi(\rho)^{\circ} \otimes_{\OO} k'$
is also cuspidal, hence isomorphic to $\pibar^{\gen} \otimes_k k'$.
\item $\pi(\rho)$ is a twist of the Steinberg representation by a character that is trivial modulo $p$.
The reduction mod $p$ of $\pi(\rho)$ then has two Jordan-H\"older constituents, isomorphic to
$\pibar^{\gen}$ and $\abs \circ \det$.  In particular $\pi(\rho)^{\circ} \otimes_{\OO} k'$
is a nonsplit extension of $(\abs \circ \det) \otimes_k k'$ by $\pibar^{\gen} \otimes_k k'$,
and is thus isomorphic to $W^{\vee} \otimes (\abs \circ \det) \otimes k'$.
\item $\pi(\rho)$ is a twist of the Steinberg representation by a character that is congruent
to $\abs \circ \det$ modulo $p$.  In this case $\pi(\rho)^{\circ} \otimes_{\OO} k'$
is isomorphic to $W^{\vee} \otimes k'$.
\item $\pi(\rho)$ is a parabolic induction, that contains a lattice whose reduction modulo $p$
is $i_B^G \abs^{\frac{1}{2}} \otimes \abs^{-{\frac{1}{2}}}$.  The reduction of
$\pi(\rho)^{\circ}$ thus has length $3$ and embeds in $V \otimes_k k'$, and is thus isomorphic to $V \otimes_k k'$.
\end{enumerate}

It is now easy to establish the modified mod $p$ local Langlands correspondence for
these $\rhobar$:

\begin{thm} \label{thm:main1}
Let $\rhobar$ be a representation of $G_F$ such that $\rhobar^{\sss} \cong 1 \oplus \overline{\omega}$.
\begin{enumerate}
\item If $\rhobar = 1 \oplus \overline{\omega}$, then $\pibar(\rhobar)$ is isomorphic to $V$.
\item If $\rhobar$ is a nonsplit extension of $1$ by $\overline{\omega}$, then $\pibar(\rhobar)$ is isomorphic
to $W^{\vee} \otimes (\abs \circ \det)$.
\item If $\rhobar$ is a nonsplit extension of $\overline{\omega}$ by $1$, then $\pibar(\rhobar)$ is
isomorphic to $W^{\vee}$.
\end{enumerate}
\end{thm}
\begin{proof}
In case (1), $\rhobar$ has a lift $\rho$ of type (4), and then $\pi(\rho)^{\circ} \otimes_{\OO} k'$ is isomoprhic to
$V \otimes_k k'$.  On the other hand $\pibar(\rhobar)$ is contained in $V$ and $\pibar(\rhobar) \otimes_k k'$
contains $\pi(\rho)^{\circ} \otimes_{\OO} k'$.  We must thus have $\pibar(\rhobar) = V$.

In case (2), $\rhobar$ has lifts of type (1) and (2), but not (3) or (4).  Thus $\pi(\rho)^{\circ} \otimes_{\OO} k'$
is contained in $W^{\vee} \otimes (\abs \circ \det) \otimes_{\OO} k'$, and the two are sometimes equal.
We must thus have $\pibar(\rhobar) = W^{\vee} \otimes (\abs \circ \det)$.  Case (3) follow from case (2) by
twisting by $\abs \circ \det$.
\end{proof}

\section{Lattices in direct sums of characters} \label{sec:lattice}

Before we turn to the case where $q \equiv 1$ mod $p$, we need a technical result.  For this section,
let $G$ be an arbitrary locally profinite group.  Let $\OO$ be a discrete valuation ring with residue field
$k$, uniformizer $\unif$, and field of fractions $K$.

Let $\chi_1$ and $\chi_2$ be two distinct characters $G \rightarrow \OO^{\times}$ that are trivial modulo $\unif$.
We will attach a class $\sigma(\chi_1,\chi_2)$ in $H^1(G,k)$ (where $G$ acts trivially on $k$) to this pair of 
characters.
Let $a$ be the largest integer such that $\chi_1$ and $\chi_2$ are congruent modulo $a$.  Then for $g$ in $G$
we define $\sigma(\chi_1,\chi_2)_g$ to be the reduction modulo $\unif$ of the element
$\frac{1}{\unif^a}(\chi_1(g) - \chi_2(g))$.

Note that $H^1(G,k)$ is isomorphic to $\Ext^1_G(1_G,1_G)$, where $1_G$ denotes the trivial character
of $G$ with values in $k$.  There is thus a bijection between lines in $H^1(G,k)$ and nonsplit extensions
of $1_G$ by $1_G$.  This bijection can be made entirely explicit as follows: let $E$ be such an extension,
let $e_1$ span the invariant line in $E$, and complete this to a $k$-basis $\{e_1,e_2\}$ of $E$.
For any $g$, $ge_2 - e_2$ is equal to $\sigma_g e_1$ for some $\sigma_g$ in $k$; the cocycle
$g \mapsto \sigma_g$ represents a class in $H^1(G,k)$ that is nontrival because $E$ is not split.  A different
choice of basis $\{e_1,e_2\}$ changes $\sigma$ by a nonzero scalar, and this gives the desired bijection of
extensions $E$ with lines in $H^1(G,k)$.

Our goal in this section is to interpret the class $\sigma(\chi_1,\chi_2)$ in terms of this isomorphism.
Let $L$ be a free $\OO$-module of rank two, with basis $e_1$ and $e_2$.  Define an action of
$G$ on $L$ by $ge_1 = \chi_1(g) e_1$ and $ge_2 = \chi_2(g) e_2$.  Let $L'$ be a $G$-stable $\OO$-lattice
in $L \otimes K$.  Then $L'/\unif L'$ is an extension of $1_G$ by $1_G$, and we have:

\begin{prop} \label{prop:class}
Suppose $L'/\unif L'$ is nonsplit.  Then $\sigma(\chi_1,\chi_2)$ generates the line in $H^1(G,k)$
corresponding to the extension $L'/\unif L'$.
\end{prop}
\begin{proof}
Since replacing $L'$ with $\unif L'$ does not change the extension $L'/\unif L'$, we may assume without
loss of generality that $L \subset L'$ but $L \not \subseteq \unif L'$.  Then the map
$$L/\unif L \rightarrow L'/\unif L'$$
has one-dimensional image.  Swapping $e_1$ and $e_2$ (and thus $\chi_1$ and $\chi_2$) if necessary
we may assume that $e_1$ generates the image of $L/\unif L$ in $L'/\unif L'$.  (Note that this only
changes $\sigma(\chi_1,\chi_2)$ by a sign.

Since $e_1$ is nonzero in $L'/\unif L'$ we may complete it to a basis $e_1, e_3$ of $L'$.
Let $b$ be the smallest integer greater than zero such that $\unif^b e_3$ lies in $L$,
and write $\unif^b e_3 = \alpha e_1 + \beta e_2$ for $\alpha, \beta \in \OO$.  We then have
$g e_3 = \chi_2(g) e_3 + \frac{1}{\unif^b} \alpha(\chi_1(g) - \chi_2(g)) e_1$.  Note that
by assumption the coefficient of $e_1$ lies in $\OO$, as $L'$ is $G$-stable.

Let $\overline{e}_1$, $\overline{e}_3$ be the images of $e_1$ and $e_2$ in $L'/\unif L'$.  The action
of $G$ fixes $\overline{e}_1$, whereas $g\overline{e}_3 = \overline{e}_3 + \sigma_g \overline{e}_1$,
where $\sigma_g$ is the reduction modulo $\unif$ of $\frac{1}{\unif b} \alpha (\chi_1(g) - \chi_2(g))$.
As $L'/\unif L'$ is nonsplit, $\sigma_g$ is nonzero for some $g$, and thus $\frac{1}{\unif^b} \alpha$
must lie in $\frac{1}{\unif^a} \OO^{\times}$.  Thus $\sigma$ is a scalar multiple of
$\sigma(\chi_1,\chi_2)$ as claimed.
\end{proof}

\section{$q \equiv 1$ mod $p$}

We now consider the case in which $q \equiv 1$ mod $p$.  In this case $\overline{\omega}$ is
the trivial character, and, so, up to twist, the only case it remains to consider is when
$\rhobar^{\sss}$ is the two-dimensional trivial representation of $G = \GL_2(F)$.  As above,
we begin by computing the appropriate essentially AIG envelope.

In this setting every subquotient of the essentially AIG envelope that contains $\pibar(\rhobar)$
has supercuspidal support given by two copies of the trivial character, and is thus isomorphic
to a subquotient of the parabolic induction $i_B^G 1_T$, where $B \subset G$ is the standard Borel,
$T$ is the standard torus, and $1_T$ is the trivial character of the torus over $k$.  This induction
has two Jordan-H\"older constituents: the trivial character $1_G$, and the Steinberg representation
$\St$ of $G$ over $k$.  We thus have $\pibar^{\gen} = \St$.

\begin{lemma} There is an isomorphism:
$$i_B^G 1_T \cong 1_G \oplus \St.$$
\end{lemma}
\begin{proof}
The restrictions $r_G^B 1_G$ and $r_G^B \St$ are both isomorphic to the trivial character $1_T$,
because the norm character $\abs$ is trivial.  We thus have:
$$\Hom_G(1_G, i_B^G 1_T) = \Hom_T(1_T,1_T)$$
$$\Hom_G(\St, i_B^G 1_T) = \Hom_T(1_T,1_T)$$
and the claim follows.
\end{proof}

\begin{lemma} The space $\Ext^1_G(1_G,\St)$ is two-dimensional.
\end{lemma}
\begin{proof}
Adjointess of parabolic induction and restriction gives an isomorphism:
$$\Ext^1_G(1_G, i_B^G 1_T) \cong \Ext^1_T(1_T,1_T)$$
and the latter is four dimensional.  On the other hand
$$\Ext^1_G(1_G, i_B^G 1_T) \cong \Ext^1_G(1_G,1_G) \oplus \Ext^1_G(1_G,\St).$$
One easily sees that $\Ext^1_G(1_G,1_G)$ is two-dimensional, and the result follows.
\end{proof}

Let $V$ be the ``universal extension'' of $1_G$ by $\St$, in other words the unique
extension of $1_G \oplus 1_G$ by $\St$ that contains every isomorphism class of
extension of $1_G$ by $\St$ (more prosaically $V$ may be constructed as the pushout:
$$
\begin{array}{ccc}
W & \rightarrow & V\\
\uparrow & & \uparrow\\
\St & \rightarrow & W'
\end{array}
$$
where $W$ and $W'$ are any two nonisomorphic extensions of $1_G$ by $\St$.)  We then have:

\begin{prop} The representation $V$ is an essentially AIG envelope of $\St$.
\end{prop}
\begin{proof}
Suppose not.  Then (just as in the $q \equiv -1$ mod $p$ case), there is an essentially AIG
representation $V'$ containing $V$ with $V'/V$ isomorphic to $1_G$.  The quotient $V'/\St$
is an extension of $1_G$ by $1_G \oplus 1_G$ on which the center of $G$ acts trivially; since $p$
is odd such an extension is split.  Thus $V'$ is an extension of $1_G \oplus 1_G \oplus 1_G$ by
$\St$; since $\Ext^1_G(1_G,\St)$ is only two dimensional we must have that $1_G$ is a direct summand of
$V'$, contradicting the fact that $V'$ is essentially AIG.
\end{proof}

It will be useful to be able to classify the nonsplit extensions of $1_G$ by $\St$.  Observe:

\begin{lemma} Let $W$ be a nonsplit extension of $1_G$ by $\St$.  Then $r_G^B W$ is a nonsplit
extension of $1_T$ by $1_T$.  (Equivalently, the map 
$$\Ext^1_G(1_G,\St) \rightarrow \Ext^1_T(1_T,1_T)$$
induced by $r_G^B$ is injective.)
\end{lemma}
\begin{proof}
Suppose $r_G^B W$ is split.  Then $\Hom_T(r_B^G W, 1_T)$ is two dimensional, so $\Hom_G(W,i_B^G T)$
is also two dimensional.  It follows that there is a surjection of $W$ onto $\St$, implying that $W$ must also
be split.
\end{proof}

It is not difficult to characterise the image of this map:
\begin{lemma}
Let $E$ be an extension of $1_G$ by $1_G$.  Then there exists an extension $W$ of 
$1_G$ by $\St$ with $r_G^B W = E$ if, and only
if, the center $Z$ of $G$ acts trivially on $E$.
\end{lemma}
\begin{proof}
The representation $W$ is essentially AIG and so $Z$ acts on $W$ by scalars, and hence trivially.
Thus the same is true of $r_B^G W$.  Thus the image of the map
$$r_G^B: \Ext^1_G(1_G,\St) \rightarrow \Ext^1_T(1_T,1_T)$$
is contained in the subspace $\Ext^1_{T/Z}(1_{T/Z},1_{T/Z})$ of $\Ext^1_T(1_T,1_T)$.  This
subspace is two-dimensional, as is the image of $r_G^B$, proving the claim.
\end{proof}

The sequence of isomorphisms: $W_F^{\mbox{\rm \tiny ab}} \cong F^{\times} \cong T/Z$ (where the last isomorphism
sends $x \in F^{\times}$ to the class of the diagonal matrix with entries $a$ and $1$) induces
a chain of isomorphisms:
$$\Ext^1_{G_F}(1_{G_F}, 1_{G_F}) \cong \Ext^1_{W_F}(1_{W_F},1_{W_F}) \cong \Ext^1_{F^{\times}}(1_{F^{\times}},1_{F^{\times}})
\cong \Ext^1_{T/Z}(1_{T/Z},1_{T/Z}).$$
Denote the composition of these morphisms by $\phi$.  We observe:

\begin{lemma}
Let $K$ be a finite extension of $\QQ_p$ with uniformizer $\unif$,
residue field $k'$ and ring of integers $\OO$.
Let ${\hat \chi}_1$ and ${\hat \chi}_2$ be distinct characters of $G_F$ with values in $\OO$,
whose reductions mod $\unif$ are trivial,
and let $\chi_1$ and $\chi_2$ be the corresponding characters of $F^{\times}$.  Then
the sequence of maps:
$$H^1(G_F, 1_G) \cong \Ext^1_{G_F}(1_G,1_G) \stackrel{\phi}{\rightarrow} \Ext^1_{T/Z}(1_{T/Z},1_{T/Z})
\cong H^1(T/Z, 1_{T/Z})$$
takes $\sigma(\chi_1,\chi_2)$ to a nonzero
multiple of the class $\sigma(\chi_1 \otimes \chi_2,\chi_2 \otimes \chi_1)$.
\end{lemma}
\begin{proof}
This is an easy computation.
\end{proof}

We are now in a position to descrive $\pibar(\rhobar)$ for each $\rhobar$.  We first enumerate
the possible lifts $\rho$ of $\rhobar$.  There are four cases:
\begin{enumerate}
\item $\rho$ is a two-dimensional representation of $G_F$ on which $G_F$ acts via a single
character ${\hat \chi}$.  In this case $\rhobar$ must be trivial.
\item $\rho$ is the direct sum of two distinct characters ${\hat \chi}_1$ and ${\hat \chi}_2$ whose reductions
are trivial.  In this case $\rhobar$ is either trivial or the unique nonsplit extension of
$1_{G_F}$ by $1_{G_F}$ of class $\sigma({\hat \chi}_1,{\hat \chi}_2)$.
\item $\rho$ is a nonsplit extension of a character $\hat \chi$ by the character $\omega \hat \chi$.
\item $\rho$ is a twist of the unique unramified extension of the trivial representation of $G_F$ over $K$
by itself.  In this case $\rhobar$ is either trivial or the unique unramified extension of $1_{G_F}$ by
$1_{G_F}$.
\end{enumerate}

The next step is to describe $\pi(\rho)^{\circ} \otimes_{\OO} k'$.  We first observe:

\begin{lemma}
Let $K$ be a finite extension of $\QQ_p$, and let $\pi$ be the irreducible
parabolic induction $i_B^G 1_{T,K}$, where $1_{T,K}$ is the one-dimensional trivial
representation of $T$ over $K$.  Then $r_G^B \pi$ is the unique nonsplit extension
of $1_{T,K}$ by $1_{T,K}$ on which the action of $T$ factors through
the quotient: $T \rightarrow T/Z \cong F^{\times} \rightarrow \ZZ$.
\end{lemma}
\begin{proof}
It is clear that the extension $r_G^B \pi$ is nonsplit, as we have isomorphisms:
$$K \cong \End_G(\pi) \cong \Hom_T(r_G^B \pi, 1_{T,K}).$$
It is also clear that the action of $T$ on $r_G^B \pi$ factors through $T/Z$.
On the other hand the representation $\pi$ is an irreducible representation
with an Iwahori fixed vector, and it is well-known that for such $\pi$,
the subgroup $\OO^{\times} \times \OO^{\times}$ of $T$ acts trivially
on $r_G^B \pi$.
%%%%if pressed, give a G-cover argument
\end{proof} 

We can now describe $\pi(\rho)^{\circ} \otimes_{\OO} k'$ in each of the above cases.
\begin{enumerate}
\item In this case $\pi(\rho)$ is a twist of $i_B^G 1_{T,K}$, and the lemma above then implies
that $\pi(\rho)^{\circ} \otimes_{\OO} k'$ is the unique extension $W$ of $1_G$ by $\St$ such that
the action of $T$ on $r_G^B W$ factors through $T \rightarrow T/Z \cong F^{\times} \rightarrow \ZZ.$
\item In this case $\pi(\rho)$ is the parabolic induction $i_B^G \chi_1 \otimes \chi_2$, where $\chi_1$
and $\chi_2$ are the characters of $F^{\times}$ arising from ${\hat \chi}_1$ and ${\hat \chi}_2$ by
local class field theory.  It follows
that $r_G^B \pi(\rho)$ is the direct sum of the characters $\chi_1 \otimes \chi_2$ and $\chi_2 \otimes \chi_1$,
and $r_G^B \pi(\rho)^{\circ} \otimes_{\OO} k'$ is then the nonsplit extension of $1_T$ by $1_T$
of class $\sigma(\chi_1 \otimes \chi_2,\chi_2 \otimes \chi_1)$.
\item In this case $\pi(\rho)$ is a twist of the Steinberg representation, and $\pi(\rho)^{\circ} \otimes_{\OO} k'$
is isomorphic to $\St \otimes_k k'$.
\item In this case $\pi(\rho)$ is a twist of $i_B^G 1_{T,K}$, and the same discussion as in case (1) applies.
\end{enumerate}

\begin{thm} \label{thm:main2}
If $\rhobar$ is trivial, then $\pibar(\rhobar) = V$.  On the other hand, if $\rhobar$ is
the nonsplit extension of $1_{G_F}$ by $1_{G_F}$ represented by $\sigma \in \Ext^1_{G_F}(1_{G_F},1_{G_F})$,
then $\pibar(\rhobar)$ is the unique nonsplit extension of $1_G$ by $\St$ such that
$r_G^B \pibar(\rhobar)$ represents the class $\phi(\sigma)$ in $\Ext^1_T(1_T,1_T)$.
\end{thm}
\begin{proof}
If $\rhobar$ is trivial, then $\rhobar$ has lifts of type (2) above for an arbitrary choice of
${\hat \chi}_1$ and ${\hat \chi}_2$.  Thus $\pi(\rho)^{\circ} \otimes_{\OO} k'$ can be an
arbitrary nonsplit extension of $1_G$ by $\St$.  As $\pibar(\rhobar) \otimes_k k'$ must contain
all of these extensions, and is contained in $V$, we must have $\pibar(\rhobar) = V$.

If $\rhobar$ is nontrivial and ramified, then $\rhobar$ has lifts of type (2) and possibly (3), but not (1)
or (4).  If $\rho$ is a lift of type (3) then $\pi(\rho)^{\circ} \otimes_{\OO} k$ is isomorphic to
$\St$ and thus tells us nothing about $\pibar(\rhobar)$.  On the other hand, if $\rho = {\hat \chi}_1 \oplus {\hat \chi}_2$
is a lift of type (2), we have $\sigma = \sigma({\hat \chi}_1,{\hat \chi}_2)$.  Then
$\pi(\rho)^{\circ} \otimes_{\OO} k'$ is the extension of $1_G$ by $\St$ corresponding to the class
$\sigma(\chi_1 \otimes \chi_2,\chi_2 \otimes \chi_1)$ in $\Ext^1_T(1_T,1_T)$.  This class is a nonzero
multiple of
$\phi(\sigma)$.  It is thus clear that $\pibar(\rhobar)$ is the extension corresponding to $\phi(\sigma)$
as claimed.

If $\rhobar$ is nontrivial but unramified, the discussion of the previous paragraph applies
but one must also consider lifts of type (4).  It suffices to check that these produce the same
extension of $1_G$ by $\St$ as the lifts of type (2); that is, that when
$\sigma$ is the class attached to an unramified nonsplit extension of $1_{G_F}$ by $1_{G_F}$,
then $\phi(\sigma)$ corresponds to the extension of $1_T$ by $1_T$ on which the action of $T$
factors through $T \rightarrow T/Z \cong F^{\times} \rightarrow \ZZ$.  This is a straightforward
calculation.
\end{proof}

\end{document}